\documentclass{article}

\title{ Gotthold Eisenstein and Philosopher John}
\author{ Franz Lemmermeyer}
\date{}

\begin{document}
\maketitle

\begin{abstract}
  Before the recent publication of the correspondence between Gauss and
  Encke, nothing was known about the role that John Taylor, a cotton
  merchant from Liverpool, had played in the life of Gotthold Eisenstein.
  In this article, we will bring together what we have discovered about
  John Taylor's life.
\end{abstract}

\section*{Eisenstein's Journey to England}

Gotthold Eisenstein
belonged, together with Dirichlet, Jacobi and Kummer,
to the generation after Gauss that shaped the theory of
numbers in the mid-19th century, and like Galois, Abel,
Riemann, Roch and Clebsch, Eisenstein died young.
Today, Eisenstein's name can be found in the Eisenstein
series, Eisenstein sums, the Eisenstein ideal, Eisenstein's
reciprocity law and in his irreducibility criterion, and he
is perhaps best known for his ingenious proofs of the
quadratic, cubic and biquadratic reciprocity laws. Eisenstein's
father Johann Konstantin Eisenstein emigrated to
England in 1840; Eisenstein and his mother followed in
June 1842, although Eisenstein's few remarks on this episode
in his autobiography \cite{EisW} belie the dramatic events
that he experienced in England. On their journey to England,
the Eisensteins passed through Hamburg; during
the Great Fire in May 1842 about a third of the houses
in the Altstadt had burned down. What we learn from
Eisenstein's account is that he was impressed by the sight
of railroad lines running right under the foundations of
houses (in London?) and by the Menai suspension bridge
in Wales: Eisenstein mentions that he undertook six sea
voyages, and that on one of them they sailed under the
tremendous suspension bridge in Anglesey, which was so
high that the Berlin Palace would easily have fitted under
its main arch.

Eisenstein also writes that he later made the acquaintance
of William Hamilton in Dublin (who would discover
the quaternions in October 1843), as well as that of the
mayor of Dublin, Daniel O'Connell (who died in 1847 in
Genoa, where he wanted to organise help for the Irish
during the Great Famine).

Eisenstein does not mention what happened after his
arrival in England; the correspondence between Encke
and Gauss reveals that he fell seriously ill and was saved
by the assistance of a certain Mr. Taylor from Liverpool.
Encke first mentions this story in his letter from 10 June
1844 (see \cite[S. 1141]{Wittmann}):


\begin{quote}
  {\em I have the honour, most esteemed Privy Counsellor,
to send you the second volume of the Berlin Observations,
and I have enclosed a collection of papers by a
young local mathematician Eisenstein, whose assessment
seems to me to be of importance in the interest
of science, and for which no better judge could ever
be found than you, the true head of this part of higher
mathematics. I would not have concerned myself with
this matter, but would have left it to Prof. Dirichlet, if
the latter had not already been away for a year.
Young Eisenstein went to the same Gymnasium
as my sons, where he did not excel at other subjects,
but surprised his teacher in mathematics\footnote{This might have
    been Karl Heinrich Schellbach.} with a curious
derivation of a series first developed by Lagrange,
without having known this paper. He was then introduced
to Prof. Dirichlet, who believed him to be an
extraordinary genius. Owing to his family being in
special financial circumstances, he followed his parents
to Liverpool two years ago, where illness and family
relations (his father seems to have speculated unsuccessfully)
reduced him to such a harassed state that
his local acquaintances thought it necessary to send
him back home. Dirichlet told me of these circumstances
and asked me to try to aid his return with the
help of my English friends. Although I did not have
any acquaintances in Liverpool, it turned out that an
important cotton merchant, Taylor, had submitted a
book on ancient Roman festivities and chronology
to the Academy, and I learned that he published articles
in journals on our comet and in particular on the
short-period comet\footnote{Encke had computed the elements of
  the comet now named after him; it has an orbital period of
  only 3.3 years.}. Based on this, I turned to this gentleman,
who was otherwise completely unknown to me, and asked for his
assistance. With a generosity that can perhaps only be found
in Englishmen, he sought out young Eisenstein, arranged for
a doctor and medicine (Eisenstein was suffering from typhoid fever),
provided him with the means for travelling to Dublin (where
Hamilton received him very obligingly) and also gave
him and his mother the necessary funds to return to
Berlin. Here, Eisenstein occupied himself with investigations
of the kind you will find among the enclosed
papers and has obtained, since other avenues were
unsuccessful, a yearly sum of money from the King
(through the efforts of Mr. v. Humboldt, as he told me),
with which this young man is completely content (in
his own words). He would like to visit G\"ottingen in the
near future in order to present himself to you and, if the
circumstances allow, profit from your teaching.
Because of his young age (he cannot yet be or must
only just be 20) and his talent, which he certainly must
possess, although I admittedly cannot tell whether it is
as great as Dirichlet's remarks would have one believe,
your judgement would be of such importance for his
future position that I entreat you with some urgency
not to deny my request. He seems to me to be pleasant
company, and his experiences do not appear to have
weighed him down, but rather shown him that he must
pull himself together, and I believe I may hope that he
will not be a burden to you. The manner in which I am
involved in these matters is not entirely comfortable to
me, since I feel indebted to Mr. Taylor and must see
how I can thank him.}
  \end{quote}
Gauss answers Encke's letter on 23 June and writes:
\begin{quote}
  {\em I would be pleased to make the acquaintance of such a promising
    young man, and I would be delighted [\ldots] if he would spend some
    time in G[\"ottingen]. It would be my pleasure if I could be of any
    assistance to him, if not by actual teaching, since he clearly has
    \underline{by far} surpassed this stage.}
\end{quote}
Gauss's judgement on Eisenstein's skills was based on the articles that
Eisenstein published in Crelle's Journal in early 1844. Gauss had only
studied one of these articles in detail, namely Eisenstein's proof of
the quadratic reciprocity law using multiple Gauss sums (see
\cite[vol. I, pp.~100--107]{EisW}), but this was sufficient to convince
Gauss of Eisentein's talent.

On 13 August 1844, Encke writes to Gauss, asking him to reply to a
letter from Taylor and reminding him what Taylor has done for Eisenstein:
\begin{quote}
  {\em Young Eisenstein, in whom you have shown so much interest, went
    to Liverpool in early 1843 (unless I'm mistaken) with his (rather
    worthless) father; there his family got into a serious plight, and
    in addition young Eisenstein became gravely ill with typhoid fever.
    Prof. Dirichlet told me about it and asked if I happened to have an
    acquaintance in Liverpool who could take an interest in the young
    man, for otherwise he would perish.

    Although I did not know Mr. Taylor (I had only heard from Prof.
    Mitzscherlich\footnote{Eilhard Mitzscherlich (1799--1863) is a
      chemist and mineralogist from Berlin.} that he was a wealthy
    cotton merchant and an amateur astronomer who published
    news about comets and was familiar with my name), I took the chance
    of approaching him about the matter. He immediately did a lot more
    than I had hoped, sending a doctor to young Eisenstein, providing
    for him by a subscription, supplying him with money for travelling
    to Dublin in order to meet Hamilton, and facilitating his return to
    Berlin.}    
\end{quote}

\newpage

\section*{John Taylor}
Who was this Mr. Taylor, cotton merchant from Liverpool,
who saved Eisenstein's life and paid for his travels
to Ireland and back to Germany? Encke writes in a letter
to Gauss from 15 August 1846:
\begin{quote}
  {\em I must confess that I think very highly of his actions, for which
    he received no compensation, and so I am sorry that he got into a
    bitter dispute with Sheepshanks concerning the building of the
    observatory, in which he (probably deservedly) drew the short straw,
    since he does not know much about modern astronomy. He knows more
    about ancient astronomy, since he has translated the first four
    books of Ovid's fastis and published them with explanations about
    the knowledge of the skies at that time.
    
    Should you thus feel inclined to answer his letter, I would sincerely
    request that you make friendly mention of his truly noble behaviour
    towards young Eisenstein, which he would value highly.}
\end{quote}
The information that we can glean from Encke's letter suffices for
identifying Taylor as the cotton merchant John Taylor from Liverpool:
both his {\em Poems and Translations} including the English translation
of Ovid's {\em fastis} \cite{Taylor} as well as his heated exchanges with
Reverend Sheepshanks in \cite{Sheepshanks} can be found online.

Richard Sheepshanks (1794--1855) was a Fellow of the Royal Society
in London since 1830. The list of persons he quarreled with is long;
his dispute with James South, in which later Charles Babbage and his
difference engine got drawn into, is described in Hoskin's article
{\em Astronomers at War} \cite{Hoskin}.

The dispute with Taylor was about the best position for the future
obervatory in Liverpool. In a letter to the Liverpool Mercury,
Taylor heaps scorn upon Sheepshanks' choice:
\begin{quote}
  {\em according to the writers of this Royal Astronomical Report, the
    proper situation for our Observatory is in the lowest point of land
    that can be found, surrounded by hills that cut off the true horizon,
    and where, in fact, there is no horizon at all, although, to be sure,
    when the trough of the river is filled with smoke and fog, which is
    commonly the case, there may seem to be a horizon at a few yards
    distance [\ldots]  where nor sun, nor moon, nor star is ever seen, or
    was ever seen, to rise or set, and where no meridian line can ever be
    drawn or determined [\ldots]}
\end{quote}
Sheepshanks prevailed, but 20 years later the observatory got closed
due to an extension of the harbor, and in 1868 the new observatory
was erected on Bidston Hill, the position that had been suggested by
Taylor.

Another source of information on John Taylor comes from the diary
of the American astronomer Maria Mitchell. In 1857 she traveled
through Europe; on August 3, 1857 (see \cite[S. 86]{Mitchell})
she delivered a letter to John Taylor and observed that he must
have been around 80:
\begin{quote}
  {\em I brought a letter from Professor Silliman to Mr. John Taylor,
    cotton merchant and astronomer; and to-day I have taken tea with him.
    He is an old man, nearly eighty I should think, but full of life, and
    talks by the hour on heathen mythology. He was the principal agent in
    the establishment of the Liverpool Observatory, but disclaims the
    honor, because it was established on so small a scale, compared with
    his own gigantic plan. Mr. Taylor has invented a little machine, for
    showing the approximate position of a comet, having the elements.
    [\ldots] He struck me as being a man of taste, but of no great profundity.}
\end{quote}
If Taylor was about 80 years old in 1857, then he must have been born
around 1777. His cometarium was studied by Beech \cite{Beech}. The
information that Taylor died  in 1857 can be found in
\cite[S.~191]{Ellison}:
\begin{quote}
  {\em One of Mr. William Ewart's strongest political supporters
    during the stirring times which preceded the passage of the first
    Reform Bill (of which Mr. Ewart became an energetic advocate)
    was Mr. John Taylor, who, along with his brother Richard, had
    commenced business as cotton broker in 1821, but who since 1826
    had been by himself. He was commonly called ``Philosopher
    John'', for besides being an active politician, both as writer
    and speaker, he was also noted as a poet and as an astronomer. He was
    the first to propose the erection of an Observatory in Liverpool,
    and out of this suggestion originated the present establishment at
    Bidston Hill. He was one of the original members of the Cotton
    Brokers' Association, and  continued in business until his death,
    which occurred in 1857.}
\end{quote}
William Ewart (1798--1869) was a liberal politician from Liverpool
fighting for the abolition of capital punishment; he voted for the
legalization of the metric system in England in 1864.

\bibliographystyle{plain}

\subsection*{Acknowledgements}
I thank the staff at the Journal of the EMS and the language editor
for considerable help with the manuscript.

\vskip 1cm

\noindent 
Franz Lemmermeyer \\
M\"orikeweg 1 \\
73489 Jagstzell \\
e-mail: hb3@ix.urz.uni-heidelberg.de 

\end{document}